\documentstyle[11pt]{amsart}

\def\opn#1#2{\def#1{\operatorname{#2}}} 
\opn\chara{char} \opn\length{\ell}
\opn\projdim{proj\,dim}
\opn\injdim{inj\,dim}
\opn\rank{rank}
\opn\depth{depth}
\opn\grade{grade}
\opn\height{height}
\opn\embdim{emb\,dim}
\opn\codim{codim}

\opn\Tr{Tr}
\opn\bigrank{big\,rank}
\opn\superheight{superheight}\opn\lcm{lcm}
\opn\trdeg{tr\,deg}%
\opn\reg{reg}
\opn\lreg{lreg}
\opn\div{div}
\opn\Div{Div}
\opn\cl{cl}
\opn\Cl{Cl}
\opn\Spec{Spec}
\opn\Supp{Supp}
\opn\supp{supp}
\opn\Sing{Sing}
\opn\Ass{Ass}
\opn\Ann{Ann}
\opn\Rad{Rad}
\opn\Soc{Soc}
\opn\Ker{Ker}
\opn\Coker{Coker}
\opn\Im{Im}
\opn\Hom{Hom}
\opn\Tor{Tor}
\opn\Ext{Ext}
\opn\End{End}
\opn\Aut{Aut}
\opn\id{id}

\opn\nat{nat}
\opn\pff{pf}
\opn\Pf{Pf}
\opn\GL{GL}
\opn\SL{SL}
\opn\mod{mod}
\opn\ord{ord}
\opn\aff{aff}
\opn\con{conv}
\opn\relint{relint}
\opn\st{st}
\opn\lk{lk}
\opn\cn{cn}
\opn\core{core}
\opn\vol{vol}
\opn\gr{gr}

\def\pot#1#2{#1[\kern-0.28ex[#2]\kern-0.28ex]}

\opn\dirlim{\underrightarrow{\lim}}
\opn\invlim{\underleftarrow{\lim}}
\def\Implies{\ifmmode\Longrightarrow \else
     \unskip${}\Longrightarrow{}$\ignorespaces\fi}
\def\implies{\ifmmode\Rightarrow \else
     \unskip${}\Rightarrow{}$\ignorespaces\fi}
\def\iff{\ifmmode\Longleftrightarrow \else
     \unskip${}\Longleftrightarrow{}$\ignorespaces\fi}

\let\:=\colon
\newtheorem{Theorem}{Theorem}[section]
\newtheorem{Lemma}[Theorem]{Lemma}

\let\epsilon=\varepsilon
\let\phi=\varphi
\let\kappa=\varkappa
\def\qed{\ifhmode\textqed\fi
   \ifmmode\ifinner\quad\qedsymbol\else\dispqed\fi\fi}
\def\textqed{\unskip\nobreak\penalty50
    \hskip2em\hbox{}\nobreak\hfil\qedsymbol
    \parfillskip=0pt \finalhyphendemerits=0}
\def\dispqed{\rlap{\qquad\qedsymbol}}

\opn\ini{in}
\opn\inm{inm}
\opn\Sym{Sym}

\begin{document}

\title{A Matrix Approach to the Rational Invariants of Certain
Classical Groups over Finite Fields of Characteristic Two}
\author{Zhongming Tang and Zhe-xian Wan}

\begin{abstract}
Let ${\mathbb F}_q$ be a finite field of characteristic two and
${\mathbb F}_q(X_1,\ldots,X_n)$ a rational function field. We use
matrix methods to obtain explicit transcendental bases of the
invariant subfields of orthogonal groups and pseudo-symplectic
groups on ${\mathbb F}_q(X_1,\ldots,X_n)$  over ${\mathbb F}_q$.
\end{abstract}

\thanks{
\indent Both authors are supported by the National Natural Science
Foundation of China}
\date{}
\maketitle

\section{Introduction}
Let ${\mathbb F}_q$ be a finite field, $n\geq 1$ an integer and $
\textmd{GL}_n({\mathbb F}_q)$ the general linear group. For any
$T\in\textmd{GL}_n({\mathbb F}_q)$, $T$ induces an ${\mathbb
F}_q$-linear action $\sigma_T$ on the rational function field
${\mathbb F}_q(X_1,\ldots,X_n)$ defined by $\sigma_T(
f(X_1,\ldots,X_n))=f(\sigma_T(X_1),\ldots,\sigma_T(X_n))$ for
all\\
$f(X_1,\ldots,X_n)\in {\mathbb F}_q(X_1,\ldots,X_n)$, where
$$
(\sigma_T(X_1),\ldots,\sigma_T(X_n))= (X_1,\ldots,X_n)T.
$$
Let $G$ be a subgroup of $\textmd{GL}_n({\mathbb F}_q)$. The
invariant subfield of $G$ on ${\mathbb F}_q(X_1,\ldots,X_n)$ is
$$
{\mathbb F}_q(X_1,\ldots,X_n)^G=\{\,f\in {\mathbb
F}_q(X_1,\ldots,X_n):\sigma_T(f)=f \mbox{ for all } T\in G\,\}.
$$
Then ${\mathbb F}_q(X_1,\ldots,X_n)$ is a finite Galois extension
of ${\mathbb F}_q(X_1,\ldots,X_n)^G$. One asks when ${\mathbb F}_q
(X_1,\ldots,X_n)^G$ is purely transcendental over ${\mathbb F}_q$
for a classical group $G$.

When $G=\textmd{GL}_n({\mathbb F}_q)$, Dickson\cite{D} gave an
affirmative answer:
$$
{\mathbb F}_q(X_1,\ldots,X_n)^{\textmd{GL}_n({\mathbb
F}_q)}={\mathbb F}_q(C_{n0},C_{n1},\ldots,C_{n,n-1}),
$$
where $C_{ni}=\frac{D_{ni}}{D_{nn}}$, $i=0,1,\ldots,n-1$, and
$$
D_{ni}=\det\left(
\begin{array}{cccccc}
X_1&X_1^q&\ldots&\widehat{X_1^{q^i}}&\ldots&X_1^{q^n}\\
X_2&X_2^q&\ldots&\widehat{X_2^{q^i}}&\ldots&X_2^{q^n}\\
&\ldots&&&\ldots&\\
X_n&X_n^q&\ldots&\widehat{X_n^{q^i}}&\ldots&X_n^{q^n}
\end{array}
\right),\,\, i=0,1,\ldots,n.
$$

Relative recently, Chu\cite{C1} considered rational invariants of
orthogonal groups and obtained similar results for $n=2,3$.
Afterwards, Cohen\cite{Co} showed the result is true when $n=4$,
and Carlisle and Kropholler\cite{CK} solved the general case. But
they all assumed that the characteristic of ${\mathbb F}_q$ is
odd. The case of characteristic two was settled by Rajaei\cite{R}
using quadratic form language.

Suppose that the characteristic of ${\mathbb F}_q$ is two. In this
paper, we will follow the route of Carlisle and
Kropholler\cite{CK} and Chu\cite{C2} to use matrix methods to give
explicit transcendental bases of the invariant subfields of
orthogonal groups on ${\mathbb F}_q(X_1,\ldots,X_n)$ over
${\mathbb F}_q$. Our transcendental bases are different from that
of \cite{R} in some cases. Moreover, we also obtain explicit
transcendental bases of the invariant subfields of
pseudo-symplectic groups on ${\mathbb F}_q(X_1,\ldots,X_n)$ over
${\mathbb F}_q$.

Following Chu\cite{C2}, some results in \cite{CK} and \cite{C2}
can be described as follows. Let $A\in\textmd{GL}_n({\mathbb
F}_q)$ be a symmetric matrix and
$$
G_A=\{T\in\textmd{GL}_n({\mathbb F}_q):\,TA{^tT}=A\},
$$
where $^tT$ denotes the transpose of $T$. Define
$$
P_{nk}=(X_1,\ldots,X_n)A\left(
\begin{array}{c}
X_1^{q^k}\\
\vdots\\
X_n^{q^k}
\end{array}
\right),\,\,k=0,1,2,\ldots.
$$
Then
$$
P_{nk}\in {\mathbb F}_q(X_1,\ldots,X_n)^{G_A},\, k=0,1,2,\ldots,
$$
and
$$
{\mathbb F}_q(X_1,\ldots,X_n)^{G_A}={\mathbb
F}_q(P_{n0},P_{n1},\ldots,P_{nn}),
$$
( cf., the proof of \cite[Theorem]{C2}). Here there is no
restriction on the characteristic of ${\mathbb F}_q$. When the
characteristic of ${\mathbb F}_q$ is odd,
$$
{\mathbb F}_q(X_1,\ldots,X_n)^{G_A}={\mathbb
F}_q(P_{n0},P_{n1},\ldots,P_{n,n-1}),
$$
i.e., $P_{nn}\in{\mathbb F}_q(P_{n0},P_{n1},\ldots,P_{n,n-1})$.
The restriction on the characteristic of ${\mathbb F}_q$ is
crucial in their proof. Similar problems for rational invariants
of unitary groups and symplectic groups have been also solved (
cf., \cite{CK} and \cite{C2}), where there is no restriction on
the characteristic of ${\mathbb F}_q$. For example, for the
symplectic group $Sp_{2\nu}({\mathbb F}_q)$ of degree $2\nu$ over
${\mathbb F}_q$ of any characteristic, we have
$$
{\mathbb F}_q(X_1,\ldots,X_n)^{Sp_{2\nu}({\mathbb F}_q)}={\mathbb
F}_q(P_{n1},P_{n2},\ldots,P_{nn}).
$$

Let $A=(a_{ij})$ be an $n\times n$ matrix over ${\mathbb F}_q$.
$A$ is said to be alternate if $a_{ij}=-a_{ji}$ and $a_{ii}=0$.
Let $A$ and $B$ be two $n\times n$ matrices, we write $A\equiv B$
if $A-B$ is alternate. The identity matrix of rank $n$ will be
denoted by $I^{(n)}$.

Suppose that the characteristic of ${\mathbb F}_q$ is two. Then,
up to isomorphisms, orthogonal groups over ${\mathbb F}_q$ are
just the following three types:
$$
O_n({\mathbb F}_q,G)=\{\,T\in\textmd{GL}_n({\mathbb F}_q):
TG{^tT}\equiv G\,\},
$$
where $G$ is one of the following
$$
\left(\begin{array}{cc} 0&I^{(\nu)}\\
&0
\end{array}\right);
\left(\begin{array}{ccc} 0&I^{(\nu)}&\\
&0&\\
&&1
\end{array}\right);
\left(\begin{array}{cccc} 0&I^{(\nu-1)}&&\\
&0&&\\
&&\alpha&1\\
&&&\alpha
\end{array}\right),
$$
with $\alpha\in{\mathbb F}_q\setminus\{x+x^2:x\in{\mathbb F}_q\}$,
which will be denoted by $G_{2\nu}$, $G_{2\nu+1}$,
$G_{2(\nu-1)+2}$, respectively. While the pseudo-symplectic group
with respect to non-alternate symmetric matrices are just the
following two types:
$$
Ps_n({\mathbb F}_q,S)=\{\,T\in\textmd{GL}_n({\mathbb
F}_q):TS{^tT}=S\,\},
$$
where $S$ is one of the following
$$
\left(
\begin{array}{ccc}
0&I^{(\nu)}&\\
I^{(\nu)}&0&\\
&&1
\end{array}
\right); \left(
\begin{array}{cccc}
0&I^{(\nu)}&&\\
I^{(\nu)}&0&&\\
&&0&1\\
&&1&1
\end{array}
\right),
$$
which will be denoted by $S_{2\nu+1}$, $S_{2\nu+2}$, respectively,
cf., Wan\cite{W}.

For the rational invariants of orthogonal groups $O_n({\mathbb
F}_q,G)$ where $G$ is $G_{2\nu}$ or $G_{2(\nu-1)+2}$, let
$A=G+{^tG}$ and $G_A$, $P_{nk}$ for $k>0$ be defined as before but
define
$$
P_{n0}=(X_1,\ldots,X_n)G\left(
\begin{array}{c}
X_1\\
\vdots\\
X_n
\end{array}
\right),
$$
(note that the former definition gives $P_{n0}=0$). Then
$O_n({\mathbb F}_q,G)$ is a subgroup of $G_A$, which is isomorphic
to $Sp_{2\nu}({\mathbb F}_q)$, and
$$
{\mathbb F}_q(X_1,\ldots,X_n)^{O_n({\mathbb F}_q,G)}={\mathbb
F}_q(P_{n0},P_{n1},\ldots,P_{nn}).
$$
The main part of the paper is devoted to show that
$P_{nn}\in{\mathbb F}_q(P_{n0},P_{n1},\ldots,P_{n,n-1})$. We need
to find a polynomial identity like $aP_{nn}=b$ where $a,b
\in{\mathbb F}_q(P_{n0},P_{n1},\ldots,P_{n,n-1})$ but $a\not=0$ as
done in odd characteristic case in \cite{CK} and \cite{C2}.
However, such an identity in odd characteristic case becomes
trivial in characteristic two ($a=0$!). To overcome this
difficulty, we first go to characteristic zero, get some identity
and cancel the factor $2$. Then, returning to characteristic two,
we obtain a required identity. In order to do this, we need some
properties of determinants which form section $2$.  The results of
rational invariants of another orthogonal group and
pseudo-symplectic groups are obtained similarly. Section 3
consists of results about orthogonal groups and section 4
discusses pseudo-symplectic groups.

\section{Preliminaries}
The following properties of determinants are useful in our
discussion. Although some of them are probably known, for the
reader's convenience, we give their proofs. The entries of
matrices will be considered as indeterminates.
\begin{Lemma}
\label{nodd} If $n$ is odd, then, in ${\mathbb Z}[a_{ij} : 1\leq
i<j\leq n]$,
$$
\det \left( \begin{array}{cccc} 0&a_{12}&\ldots&a_{1n}\\
a_{12}&0&\ldots&a_{2n}\\
&\ldots&\ldots&\\
a_{1n}&a_{2n}&\ldots&0 \end{array} \right) \equiv 0 (\mod 2).
$$
\end{Lemma}
\begin{pf} It is because that, in ${\mathbb Z}_2[a_{ij} : 1\leq
i<j\leq n]$,
$$
\det \left( \begin{array}{cccc} 0&a_{12}&\ldots&a_{1n}\\
a_{12}&0&\ldots&a_{2n}\\
&\ldots&\ldots&\\
a_{1n}&a_{2n}&\ldots&0 \end{array} \right)=0.
$$
\end{pf}

\begin{Lemma}
\label{neven} If $n$ is even, then, in ${\mathbb Z}[a_{ij} : 1\leq
i\leq j\leq n]$,
$$
\det \left( \begin{array}{cccc} 2a_{11}&a_{12}&\ldots&a_{1n}\\
a_{12}&2a_{22}&\ldots&a_{2n}\\
&\ldots&\ldots&\\
a_{1n}&a_{2n}&\ldots&2a_{nn} \end{array} \right)\equiv
\det \left( \begin{array}{cccc} 0&a_{12}&\ldots&a_{1n}\\
a_{12}&0&\ldots&a_{2n}\\
&\ldots&\ldots&\\
a_{1n}&a_{2n}&\ldots&0 \end{array} \right) (\mod 4).
$$
\end{Lemma}
\begin{pf} Since
\begin{eqnarray*}
&&
\det \left( \begin{array}{cccc} 2a_{11}&a_{12}&\ldots&a_{1n}\\
a_{12}&2a_{22}&\ldots&a_{2n}\\
&\ldots&\ldots&\\
a_{1n}&a_{2n}&\ldots&2a_{nn} \end{array} \right)\\
&=& 2a_{11}\det \left( \begin{array}{cccc} 2a_{22}&a_{23}&\ldots&a_{2n}\\
a_{23}&2a_{33}&\ldots&a_{3n}\\
&\ldots&\ldots&\\
a_{2n}&a_{3n}&\ldots&2a_{nn} \end{array} \right)
+\det \left( \begin{array}{cccc} 0&a_{12}&\ldots&a_{1n}\\
a_{12}&2a_{22}&\ldots&a_{2n}\\
&\ldots&\ldots&\\
a_{1n}&a_{2n}&\ldots&2a_{nn} \end{array} \right)\\
&\equiv&
2a_{11}\det \left( \begin{array}{cccc} 0&a_{23}&\ldots&a_{2n}\\
a_{23}&0&\ldots&a_{3n}\\
&\ldots&\ldots&\\
a_{2n}&a_{3n}&\ldots&0 \end{array} \right)+\det \left( \begin{array}{cccc} 0&a_{12}&\ldots&a_{1n}\\
a_{12}&2a_{22}&\ldots&a_{2n}\\
&\ldots&\ldots&\\
a_{1n}&a_{2n}&\ldots&2a_{nn} \end{array} \right) (\mod 4),
\end{eqnarray*}
it follows from \ref{nodd} that
$$
\det \left( \begin{array}{cccc} 2a_{11}&a_{12}&\ldots&a_{1n}\\
a_{12}&2a_{22}&\ldots&a_{2n}\\
&\ldots&\ldots&\\
a_{1n}&a_{2n}&\ldots&2a_{nn} \end{array} \right)\equiv
\det \left( \begin{array}{cccc} 0&a_{12}&\ldots&a_{1n}\\
a_{12}&2a_{22}&\ldots&a_{2n}\\
&\ldots&\ldots&\\
a_{1n}&a_{2n}&\ldots&2a_{nn} \end{array} \right) (\mod 4).
$$
Then the result follows by induction.
\end{pf}

\begin{Lemma}
\label{difference} If $n=2\nu$ is even, then, in ${\mathbb
Z}[a_{ij} : 1\leq i<j\leq n]$,
$$
\det\left( \begin{array}{ccccc} 0&a_{12}&a_{13}&\ldots&a_{1n}\\
-a_{12}&0&a_{23}&\ldots&a_{2n}\\
-a_{13}&-a_{23}&0&\ldots&a_{3n}\\
&\ldots&&\ldots&\\
-a_{1n}&-a_{2n}&-a_{3n}&\ldots&0 \end{array} \right) \equiv
(-1)^{\nu}\det
\left( \begin{array}{ccccc} 0&a_{12}&a_{13}&\ldots&a_{1n}\\
a_{12}&0&a_{23}&\ldots&a_{2n}\\
a_{13}&a_{23}&0&\ldots&a_{3n}\\
&\ldots&&\ldots&\\
a_{1n}&a_{2n}&a_{3n}&\ldots&0 \end{array} \right)(\mod 4).
$$
\end{Lemma}

\begin{pf}
By induction on $\nu$. The case $\nu=1$ is clear. Suppose that
$\nu>1$ and the result is true for $\nu-1$. Set
$$
A=\left( \begin{array}{ccccc} 0&a_{12}&a_{13}&\ldots&a_{1n}\\
-a_{12}&0&a_{23}&\ldots&a_{2n}\\
-a_{13}&-a_{23}&0&\ldots&a_{3n}\\
&\ldots&&\ldots&\\
-a_{1n}&-a_{2n}&-a_{3n}&\ldots&0 \end{array} \right),
B=\left( \begin{array}{ccccc} 0&a_{12}&a_{13}&\ldots&a_{1n}\\
a_{12}&0&a_{23}&\ldots&a_{2n}\\
a_{13}&a_{23}&0&\ldots&a_{3n}\\
&\ldots&&\ldots&\\
a_{1n}&a_{2n}&a_{3n}&\ldots&0 \end{array} \right).
$$
Let us consider $a_{ij}$ as indeterminates and work in the field
${\mathbb Q}(a_{ij}:1\leq i<j\leq n)$. Let $E$ be the matrix
obtained from the $2\nu\times 2\nu$ identity matrix by replacing
the zeros at $(2,3), (2,4), \ldots, (2,n)$ positions by
$-\frac{a_{13}}{a_{12}},-\frac{a_{14}}{a_{12}},\ldots,-\frac{a_{1n}}{a_{12}}$,
respectively. Then
$$
{^tE}AE=\left( \begin{array}{ccccc} 0&a_{12}&0&\ldots&0\\
-a_{12}&0&a_{23}&\ldots&a_{2n}\\
0&-a_{23}&&&\\
&\ldots&&A_1&\\
0&-a_{2n}&&& \end{array} \right),
{^tE}BE=\left( \begin{array}{ccccc} 0&a_{12}&0&\ldots&0\\
a_{12}&0&a_{23}&\ldots&a_{2n}\\
0&a_{23}&&&\\
&\ldots&&B_1&\\
0&a_{2n}&&& \end{array} \right),
$$
where
$$
A_1=\left( \begin{array}{cccc} 0&a'_{34}&\ldots&a'_{3n}\\
-a'_{34}&0&\ldots&a'_{4n}\\
&\ldots&\ldots&\\
-a'_{3n}&-a'_{4n}&\ldots&0 \end{array} \right),
B_1=\left( \begin{array}{cccc} a^{''}_{33}&a^{''}_{34}&\ldots&a^{''}_{3n}\\
a^{''}_{34}&a^{''}_{44}&\ldots&a^{''}_{4n}\\
&\ldots&\ldots&\\
a^{''}_{3n}&a^{''}_{4n}&\ldots&a^{''}_{nn} \end{array} \right),
$$
and
$a'_{ij}=a_{ij}-\frac{a_{1i}}{a_{12}}a_{2j}+\frac{a_{1j}}{a_{12}}a_{2i}$,
$a^{''}_{ij}=-2\frac{a_{1j}}{a_{12}}a_{2i}+a'_{ij}$ for $i<j$,
$a^{''}_{ii}=-2\frac{a_{1i}}{a_{12}}a_{2i}$. Note that
$a_{12}a'_{ij},a_{12}a^{''}_{ij}\in{\mathbb Z}[a_{ij}:1\leq
i<j\leq n]$.

Clearly $\det(A)=a^2_{12}\det(A_1)$ and
$\det(B)=-a^2_{12}\det(B_1)$. Let
$$
B_2=\left( \begin{array}{cccc} 0&a'_{34}&\ldots&a'_{3n}\\
a'_{34}&0&\ldots&a'_{4n}\\
&\ldots&\ldots&\\
a'_{3n}&a'_{4n}&\ldots&0 \end{array} \right).
$$

By induction assumption,
$a_{12}^{n-2}(\det(A_1)-(-1)^{\nu-1}\det(B_2))\in 4{\mathbb
Z}[a_{ij}:1\leq i<j\leq n]$. Note that
\begin{eqnarray*}
&&a^{n-4}_{12}(\det(A)-(-1)^{\nu}\det(B))\\
&=&a_{12}^{n-2}(\det(A_1)-(-1)^{\nu-1}\det(B_1))\\
&=&a_{12}^{n-2}(\det(A_1)-(-1)^{\nu-1}\det(B_2))+(-1)^{\nu-1}a_{12}^{n-2}(\det(B_2)-\det(B_1)).
\end{eqnarray*}
Thus it is enough to show that
$a_{12}^{n-2}(\det(B_1)-\det(B_2))\in 4{\mathbb Z}[a_{ij}:1\leq
i<j\leq n]$. Let $c_{ij}=a'_{ij}$ for $3\leq i<j\leq n$, and
$d_{ij}=-\frac{a_{1i}}{a_{12}}a_{2j}$ for $3\leq i\leq j\leq n$.
Then
\begin{eqnarray*}
B_1&=&
\left( \begin{array}{cccc} 2d_{33}&2d_{34}+c_{34}&\ldots&2d_{3n}+c_{3n}\\
2d_{34}+c_{34}&2d_{44}&\ldots&2d_{4n}+c_{4n}\\
&\ldots&\ldots&\\
2d_{3n}+c_{3n}&2d_{4n}+c_{4n}&\ldots&2d_{nn} \end{array}
\right),\\
B_2&=&\left( \begin{array}{cccc} 0&c_{34}&\ldots&c_{3n}\\
c_{34}&0&\ldots&c_{4n}\\
&\ldots&\ldots&\\
c_{3n}&c_{4n}&\ldots&0 \end{array} \right).
\end{eqnarray*}
Let us show that $\det(B_1)-\det(B_2)\in 4{\mathbb
Z}[d_{ij},c_{ij}:3\leq i\leq j\leq n],$ where $c_{ii}=0$, then the
lemma follows.

Note that, by \ref{neven},
$$
\det(B_1)-\det\left( \begin{array}{cccc} 0&2d_{34}+c_{34}&\ldots&2d_{3n}+c_{3n}\\
2d_{34}+c_{34}&0&\ldots&2d_{4n}+c_{4n}\\
&\ldots&\ldots&\\
2d_{3n}+c_{3n}&2d_{4n}+c_{4n}&\ldots&0 \end{array} \right) \in 4
{\mathbb Z}[d_{ij},c_{ij}:3\leq i\leq j\leq n],
$$
hence, we may assume that $d_{ii}=0$, $i=3,\ldots,n$. Then, in
${\mathbb Z}[d_{ij},c_{ij}:3\leq i\leq j\leq n]$,
\begin{eqnarray*}
&&\det(B_1)\\
&=&
\det\left( \begin{array}{cccc} 0&2d_{34}&\ldots&2d_{3n}\\
2d_{34}&0&\ldots&2d_{4n}+c_{4n}\\
&\ldots&\ldots&\\
2d_{3n}&2d_{4n}+c_{4n}&\ldots&0 \end{array} \right)
+\det\left( \begin{array}{cccc} 0&2d_{34}&\ldots&2d_{3n}\\
c_{34}&0&\ldots&2d_{4n}+c_{4n}\\
&\ldots&\ldots&\\
c_{3n}&2d_{4n}+c_{4n}&\ldots&0 \end{array} \right)\\
&&
+\det\left( \begin{array}{cccc} 0&c_{34}&\ldots&c_{3n}\\
2d_{34}&0&\ldots&2d_{4n}+c_{4n}\\
&\ldots&\ldots&\\
2d_{3n}&2d_{4n}+c_{4n}&\ldots&0 \end{array} \right)
+\det\left( \begin{array}{cccc} 0&c_{34}&\ldots&c_{3n}\\
c_{34}&0&\ldots&2d_{4n}+c_{4n}\\
&\ldots&\ldots&\\
c_{3n}&2d_{4n}+c_{4n}&\ldots&0 \end{array} \right)\\
&\equiv&\det\left( \begin{array}{cccc} 0&c_{34}&\ldots&c_{3n}\\
c_{34}&0&\ldots&2d_{4n}+c_{4n}\\
&\ldots&\ldots&\\
c_{3n}&2d_{4n}+c_{4n}&\ldots&0 \end{array} \right) (\mod 4),
\end{eqnarray*}
thus, continuing reductively we get that, in ${\mathbb
Z}[d_{ij},c_{ij}:3\leq i\leq j\leq n]$,
$$
\det(B_1)\equiv
\det\left( \begin{array}{cccc} 0&c_{34}&\ldots&c_{3n}\\
c_{34}&0&\ldots&c_{4n}\\
&\ldots&\ldots&\\
c_{3n}&c_{4n}&\ldots&0 \end{array} \right) (\mod 4),
$$
as required.
\end{pf}

Let $A=(a_{ij})$ be an $n\times n$ alternate matrix. When $n$ is
even, one can define the Pfaffian of $A$, denoted by $\Pf(A)$,
cf., \cite[\S 5, no.2]{B}. When $n$ is odd, set $\Pf(A)=0$. The
following properties of the Pfaffians are well-known.
\begin{itemize}
\item[(i)] $\det(A)=(\Pf(A))^2$. \item[(ii)] Interchanging two
rows and same columns simultaneously, the Pfaffian changes by
$-1$. \item[(iii)] Pfaffians can be developed along a row:
$$
\Pf(A)=\sum_{j=1}^n(-1)^{i+j-1}\sigma(i,j)a_{ij}\Pf(A_{i,j}),\,\,i=1,\ldots,n,
$$
where $\sigma(i,j)$ is the sign of $j-i$ and $A_{i,j}$ is the
matrix obtained from $A$ by deleting its $i$-th and $j$-th rows
and columns together.
\end{itemize}

We also need the following property of Pfaffians.

\begin{Lemma}
\label{multiplicity} Let $n>0$ be an even integer and $C$ an
$(n-1)\times(n-1)$ alternate matrix. Then, for any
$a_{12},\ldots,a_{1n},b_{12},\ldots,b_{1n}$,
$$
\Pf\left(\begin{array}{cc} 0&a_{12}\, \cdots\,a_{1n}\\
\begin{array}{c}
-a_{12}\\
\vdots\\
-a_{1n}
\end{array}&C
\end{array}\right)\cdot
\Pf\left(\begin{array}{cc} 0&b_{12}\, \cdots\,b_{1n}\\
\begin{array}{c}
-b_{12}\\
\vdots\\
-b_{1n}
\end{array}&C
\end{array}\right)
= \det
\left(\begin{array}{cc} 0&a_{12}\, \cdots\,a_{1n}\\
\begin{array}{c}
-b_{12}\\
\vdots\\
-b_{1n}
\end{array}&C
\end{array}\right).
$$
\end{Lemma}

\begin{pf}
Consider $a_{1i}$, $b_{1j}$ and the non-diagonal elements of $C$
as indeterminates and work in the rational field generated by them
over ${\mathbb Q}$. Let
$$
A=\left(\begin{array}{cc} 0&a_{12}\, \cdots\,a_{1n}\\
\begin{array}{c}
-a_{12}\\
\vdots\\
-a_{1n}
\end{array}&C
\end{array}\right),B=\left(\begin{array}{cc} 0&b_{12}\, \cdots\,b_{1n}\\
\begin{array}{c}
-b_{12}\\
\vdots\\
-b_{1n}
\end{array}&C
\end{array}\right),D=\left(\begin{array}{cc} 0&a_{12}\, \cdots\,a_{1n}\\
\begin{array}{c}
-b_{12}\\
\vdots\\
-b_{1n}
\end{array}&C
\end{array}\right).
$$
We first use induction on $\frac{n}{2}$ to show that
$\det(A)\cdot\det(B)=(\det(D))^2$. When $n=2$, it is clear. Now
assume that $n\geq 4$. Write
$$
C=\left(\begin{array}{cc} C_{11}&C_{12}\\
-{^tC_{12}}&C_{22}\end{array}\right),
$$
where $C_{11}$ is $(n-3)\times(n-3)$, $C_{12}$ is $(n-3)\times 2$
and $C_{22}$ is $2\times 2$, and
$$
A_1=(a_{12}\cdots a_{1,n-2}),A_2=(a_{1,n-1}\, a_{1n}),
B_1=(b_{12}\cdots b_{1,n-2}),B_2=(b_{1,n-1}\, b_{1n}).
$$
Then
\begin{eqnarray*} &&
\left(\begin{array}{ccc} 1&&-A_2C_{22}^{-1}\\
&I^{(n-3)}&-C_{12}C_{22}^{-1}\\
&&I^{(2)}
\end{array}
\right)A\left(\begin{array}{ccc} 1&&\\
&I^{(n-3)}&\\
C_{22}^{-1}{^tA_2}&C_{22}^{-1}{^tC_{12}}&I^{(2)}
\end{array}
\right)\\
&=&\left(\begin{array}{ccc}0&A_1+A_2C_{22}^{-1}{^tC_{12}}&0\\
-{^tA_1}+C_{12}C_{22}^{-1}{^tA_2}&C_{11}+C_{12}C_{22}^{-1}{^tC_{12}}&0\\
0&0&C_{22}
\end{array}\right).
\end{eqnarray*}
Thus
$$
\det(A)=\det(C_{22})\cdot\det \left(\begin{array}{cc}
0&A_1+A_2C_{22}^{-1}{^tC_{12}}\\
-{^tA_1}+C_{12}C_{22}^{-1}{^tA_2}&C_{11}+C_{12}C_{22}^{-1}{^tC_{12}}
\end{array}
\right).
$$
Similarly,
\begin{eqnarray*}
\det(B)&=&\det(C_{22})\cdot\det \left(\begin{array}{cc}
0&B_1+B_2C_{22}^{-1}{^tC_{12}}\\
-{^tB_1}+C_{12}C_{22}^{-1}{^tB_2}&C_{11}+C_{12}C_{22}^{-1}{^tC_{12}}
\end{array}
\right)\\
\det(D)&=&\det(C_{22})\cdot\det \left(\begin{array}{cc}
A_2C_{22}^{-1}{^tB_2}&A_1+A_2C_{22}^{-1}{^tC_{12}}\\
-{^tB_1}+C_{12}C_{22}^{-1}{^tB_2}&C_{11}+C_{12}C_{22}^{-1}{^tC_{12}}
\end{array}
\right).
\end{eqnarray*}
Since $C_{11}+C_{12}C_{22}^{-1}{^tC_{12}}$ is alternate of odd
order $n-3$, its determinant is $0$, hence
\begin{eqnarray*}
&& \det\left(\begin{array}{cc}
A_2C_{22}^{-1}{^tB_2}&A_1+A_2C_{22}^{-1}{^tC_{12}}\\
-{^tB_1}+C_{12}C_{22}^{-1}{^tB_2}&C_{11}+C_{12}C_{22}^{-1}{^tC_{12}}
\end{array}
\right)\\
&=&\det \left(\begin{array}{cc}
0&A_1+A_2C_{22}^{-1}{^tC_{12}}\\
-{^tB_1}+C_{12}C_{22}^{-1}{^tB_2}&C_{11}+C_{12}C_{22}^{-1}{^tC_{12}}
\end{array}
\right).
\end{eqnarray*}
By induction hypothesis
$$
\det\left(\begin{array}{cc}
0&A_1+A_2C_{22}^{-1}{^tC_{12}}\\
-{^tA_1}+C_{12}C_{22}^{-1}{^tA_2}&C_{11}+C_{12}C_{22}^{-1}{^tC_{12}}
\end{array}
\right)\cdot\det \left(\begin{array}{cc}
0&B_1+B_2C_{22}^{-1}{^tC_{12}}\\
-{^tB_1}+C_{12}C_{22}^{-1}{^tB_2}&C_{11}+C_{12}C_{22}^{-1}{^tC_{12}}
\end{array}
\right) $$ $$ =\left(\det \left(\begin{array}{cc}
0&A_1+A_2C_{22}^{-1}{^tC_{12}}\\
-{^tB_1}+C_{12}C_{22}^{-1}{^tB_2}&C_{11}+C_{12}C_{22}^{-1}{^tC_{12}}
\end{array}
\right)\right)^2.
$$
Therefore $\det(A)\cdot\det(B)=(\det(D))^2$. It follows that
$\Pf(A)\cdot\Pf(B)=\pm\det(D)$.

Now, we use one term to determine the sigh in front of
$(\det(D))^2$. Let $C=(c_{ij})$. Since the  sigh of term
$a_{12}c_{23}c_{45}\cdots c_{n-2,n-1}$ in $\Pf(A)$ is $1$, the
sign of term $b_{12}c_{23}c_{45}\cdots c_{n-2,n-1}$ in $\Pf(B)$ is
$1$ and the sigh of the corresponding term
$a_{12}b_{12}c_{23}^2c_{45}^2\cdots c_{n-2,n-1}^2$ in $\det(D)$ is
also $1$, it follows that the above sigh in front of $\det(D)$ is
$1$, i.e., $\Pf(A)\cdot\Pf(B)=\det(D)$.
\end{pf}

\section{Rational Invariants of Orthogonal Groups}
Consider first the case when $n=2\nu$ is even. We begin by
deducing some polynomial identities over ${\mathbb Z}$. Let $q$ be
a power of a prime $p$ and $f=\sum
a_{\alpha_1\cdots\alpha_n}X_1^{\alpha_1}\cdots X_n^{\alpha_n}\in
{\mathbb Z}[X_1,\ldots,X_n]$. Denote $\sum
a_{\alpha_1\cdots\alpha_n} X_1^{\alpha_1q^r}\cdots
X_n^{\alpha_nq^r}$ by $f^{(q^r)}$. For any matrix $C=(c_{ij})$
over ${\mathbb Z}[X_1,\ldots,X_n]$, let
$C^{(q^r)}=(c_{ij}^{(q^r)})$. Imbed ${\mathbb Z}_p$ in ${\mathbb
F}_q$ as a subring. Under the natural map from ${\mathbb Z}$ to
${\mathbb Z}_p$, we will use the same $f$ to denote the natural
image of $f\in {\mathbb Z}[X_1,\ldots,X_n]$ in ${\mathbb
Z}_p[X_1,\ldots,X_n]$. By our convention, $f^{(q^r)}=f^{q^r}$ in
${\mathbb Z}_p[X_1,\ldots,X_n]$ for any $f\in {\mathbb
Z}[X_1,\ldots,X_n]$. Let $G=G_{2\nu}$ or $G_{2(\nu-1)+2}$, but now
the $\alpha$ in $G_{2(\nu-1)+2}$ is regarded as an indeterminate.
Set $\alpha=0$ when $G=G_{2\nu}$. Let $A=G+{^tG}$ and
$\overline{A}=G-{^tG}$. Then $\det
(A)=(-1)^{\nu-1}4\alpha^2+(-1)^{\nu}$ and $\det(\overline{A})=1$
thus both $A$ and $\overline{A}$ are nonsingular in ${\mathbb
Z}(\alpha)$ and also in ${\mathbb Z}_p(\alpha)$. In ${\mathbb
Z}[X_1,\ldots,X_n]$, define
\begin{eqnarray}
P_{n0}&=&(X_1,\ldots,X_n)G\left(
\begin{array}{c}
X_1\\
\vdots\\
X_n\end{array}
\right),\\
P_{nk}&=&(X_1,\ldots,X_n)A\left(
\begin{array}{c}
X_1^{q^k}\\
\vdots\\
X_n^{q^k}
\end{array}
\right),\,\,k=1,2,\ldots.
\end{eqnarray}

Let
$$
B=\left(
\begin{array}{cccc}
X_1&X_2&\cdots&X_n\\
X^q_1&X^q_2&\cdots&X^q_n\\
&\cdots&\cdots&\\
X^{q^{n-1}}_1&X^{q^{n-1}}_2&\cdots&X^{q^{n-1}}_n
\end{array}\right).
$$
and
$$
P=\left( \begin{array}{ccccc} 0&P_{n1}&P_{n2}&\ldots&P_{n,n-1}\\
-P_{n1}&0&P_{n1}^{(q)}&\ldots&P_{n,n-2}^{(q)}\\
-P_{n2}&-P_{n1}^{(q)}&0&\ldots&P_{n,n-3}^{(q^2)}\\
&\ldots&&\ldots&\\
-P_{n,n-1}&-P_{n,n-2}^{(q)}&-P_{n,n-3}^{(q^2)}&\ldots&0
\end{array} \right).
$$

Then we have

\begin{Lemma}
\label{pftimes}
$$ \Pf(P)\cdot \Pf(P^{(q)})=\det
\left( \begin{array}{ccccc} 0&P_{n1}&P_{n2}&\ldots&P_{n,n-1}\\
-P_{n,n-1}^{(q)}&0&P_{n1}^{(q)}&\ldots&P_{n,n-2}^{(q)}\\
-P_{n,n-2}^{(q^2)}&-P_{n1}^{(q)}&0&\ldots&P_{n,n-3}^{(q^2)}\\
&\ldots&&\ldots&\\
-P_{n1}^{(q^{n-1})}&-P_{n,n-2}^{(q)}&-P_{n,n-3}^{(q^2)}&\ldots&0
\end{array} \right).
$$
\end{Lemma}

\begin{pf} Note that
\begin{eqnarray*}
\Pf(P^{(q)})&=& \Pf\left( \begin{array}{ccccc} 0&P_{n1}^{(q)}&P_{n2}^{(q)}&\ldots&P_{n,n-1}^{(q)}\\
-P_{n1}^{(q)}&0&P_{n1}^{(q^2)}&\ldots&P_{n,n-2}^{(q^2)}\\
-P_{n2}^{(q)}&-P_{n1}^{(q^2)}&0&\ldots&P_{n,n-3}^{(q^3)}\\
&\ldots&&\ldots&\\
-P_{n,n-1}^{(q)}&-P_{n,n-2}^{(q^2)}&-P_{n,n-3}^{(q^3)}&\ldots&0
\end{array} \right)\\
&=&
-\Pf\left( \begin{array}{ccccc} 0&-P_{n,n-1}^{(q)}&-P_{n,n-2}^{(q^2)}&\ldots&-P_{n1}^{(q^{n-1})}\\
P_{n,n-1}^{(q)}&0&P_{n1}^{(q)}&\ldots&P_{n,n-2}^{(q)}\\
P_{n,n-2}^{(q^2)}&-P_{n1}^{(q)}&0&\ldots&P_{n,n-3}^{(q^2)}\\
&\ldots&&\ldots&\\
P_{n1}^{(q^{n-1})}&-P_{n,n-2}^{(q)}&-P_{n,n-3}^{(q^2)}&\ldots&0
\end{array} \right),
\end{eqnarray*}
by property (ii) of Pfaffians. Then, our lemma follows from
\ref{multiplicity}.
\end{pf}

Let $E=(e_{ij})$ be an $m\times m$ matrix over ${\mathbb
Z}[X_1,\ldots,X_n]$. For convenience, we will use ${\mathbb Z}[E]$
to denote ${\mathbb Z}[e_{ij}:i,j=1,\ldots,m]$.

\begin{Lemma}
\label{mainlemma} There is a polynomial $S\in{\mathbb
Z}[X_1,\ldots,X_n]$ such that
$$
\det(B)\equiv2S+\Pf(P) (\mod 4)
$$
and
$$
S(S+\Pf(P))\equiv f (\mod 2)
$$
for some $f\in{\mathbb
Z}[\alpha,P_{n0},P_{n0}^{(q)},\ldots,P_{n0}^{(q^{n-1})},P]$.
\end{Lemma}

\begin{pf} Let
\begin{eqnarray*}
\overline{P}_{nk}&=&(X_1,\ldots,X_n)\overline{A}\left(
\begin{array}{c}
X_1^{q^k}\\
\vdots\\
X_n^{q^k}
\end{array}
\right),\\
P'_{nk}&=&(X_1,\ldots,X_n){^tG}\left(
\begin{array}{c}
X_1^{q^k}\\
\vdots\\
X_n^{q^k}
\end{array}
\right),\,\,k=1,2,\ldots,
\end{eqnarray*}
and
$$
\overline{P}=\left( \begin{array}{ccccc} 0&\overline{P}_{n1}&\overline{P}_{n2}&\ldots&\overline{P}_{n,n-1}\\
-\overline{P}_{n1}&0&\overline{P}_{n1}^{(q)}&\ldots&\overline{P}_{n,n-2}^{(q)}\\
-\overline{P}_{n2}&-\overline{P}_{n1}^{(q)}&0&\ldots&\overline{P}_{n,n-3}^{(q^2)}\\
&\ldots&&\ldots&\\
-\overline{P}_{n,n-1}&-\overline{P}_{n,n-2}^{(q)}&-\overline{P}_{n,n-3}^{(q^2)}&\ldots&0
\end{array} \right).
$$
Clearly, $B\overline{A}{^tB}=\overline{P}$. Since
$\det(\overline{A})=1$,
$(\det(B))^2=\det(B\overline{A}{^tB})=\det(\overline{P})=(\Pf(\overline{P}))^2$.
Thus $\det(B)=\pm\Pf(\overline{P})$. Clearly,
$\overline{P}_{nk}=P_{nk}-2P'_{nk}$. Therefore
\begin{eqnarray*}
&&
\Pf(\overline{P})\\
&=& \Pf
\left( \begin{array}{ccccc} 0&P_{n1}-2P'_{n1}&P_{n2}-2P'_{n2}&\ldots&P_{n,n-1}-2P'_{n,n-1}\\
-P_{n1}+2P'_{n1}&0&P_{n1}^{(q)}-2P^{'(q)}_{n1}&\ldots&P_{n,n-2}^{(q)}-2P^{'(q)}_{n,n-2}\\
-P_{n2}+2P'_{n2}&-P_{n1}^{(q)}+2P^{'(q)}_{n1}&0&\ldots&P_{n,n-3}^{(q^2)}-2P^{'(q^2)}_{n,n-3}\\
&\ldots&&\ldots&\\
-P_{n,n-1}+2P'_{n,n-1}&-P_{n,n-2}^{(q)}+2P^{'(q)}_{n,n-2}&-P_{n,n-3}^{(q^2)}+2P^{'(q^2)}_{n,n-3}&\ldots&0
\end{array} \right)\\
&\equiv&\Pf(P) (\mod 2),
\end{eqnarray*}
which implies
$$
\det(B)\equiv 2S+\Pf(P) (\mod 4),
$$
where $S\in {\mathbb Z}[X_1,\ldots,X_n]$, and
$$
(\det(B))^2\equiv 4S^2+4\Pf(P)S+(\Pf(P))^2 (\mod 8).
$$
On the other hand, by \ref{neven},
\begin{eqnarray*}
(\det(B))^2&=&\frac{1}{\det A}\det(BA{^tB})\\
&=&\frac{1}{\det A}\det
\left( \begin{array}{cccc} 2P_{n0}&P_{n1}&\ldots&P_{n,n-1}\\
P_{n1}&2P_{n0}^{(q)}&\ldots&P_{n,n-2}^{(q)}\\
&\ldots&\ldots&\\
P_{n,n-1}&P_{n,n-2}^{(q)}&\ldots&2P_{n0}^{(q^{n-1})}
\end{array} \right)\\
&=& \frac{1}{\det A}(4g+
\det\left( \begin{array}{cccc} 0&P_{n1}&\ldots&P_{n,n-1}\\
P_{n1}&0&\ldots&P_{n,n-2}^{(q)}\\
&\ldots&\ldots&\\
P_{n,n-1}&P_{n,n-2}^{(q)}&\ldots&0
\end{array} \right)),
\end{eqnarray*}
with some $g\in{\mathbb
Z}[P_{n0},P_{n0}^{(q)},\ldots,P_{n0}^{(q^{n-1})},P]$. Thus,
\begin{eqnarray*}
&&4g+
\det\left( \begin{array}{cccc} 0&P_{n1}&\ldots&P_{n,n-1}\\
P_{n1}&0&\ldots&P_{n,n-2}^{(q)}\\
&\ldots&\ldots&\\
P_{n,n-1}&P_{n,n-2}^{(q)}&\ldots&0
\end{array} \right)\\
&\equiv&\det(A)(4S^2+4\Pf(P)S+(\Pf(P))^2) (\mod 8)\\
&\equiv&((-1)^{\nu-1}4\alpha^2+(-1)^{\nu})(4S^2+4\Pf(P)S+(\Pf(P))^2) (\mod 8)\\
&\equiv&(-1)^{\nu}(4S^2+4\Pf(P)S+(\Pf(P))^2)+(-1)^{\nu-1}4\alpha^2(\Pf(P))^2
(\mod 8).
\end{eqnarray*}
But, by \ref{difference},
$$
(\Pf(P))^2-(-1)^{\nu}
\det\left( \begin{array}{cccc} 0&P_{n1}&\ldots&P_{n,n-1}\\
P_{n1}&0&\ldots&P_{n,n-2}^{(q)}\\
&\ldots&\ldots&\\
P_{n,n-1}&P_{n,n-2}^{(q)}&\ldots&0
\end{array} \right)\in 4{\mathbb Z}[P].
$$
Hence $ S^2+\Pf(P)S\equiv f (\mod 2)$, for some $f\in{\mathbb
Z}[\alpha,P_{n0},P_{n0}^{(q)},\ldots,P_{n0}^{(q^{n-1})},P]$, as
required.
\end{pf}

Now we can prove the main theorem.

\begin{Theorem}
\label{mainthm} Let ${\mathbb F}_q$ be a finite field of
characteristic two, $n=2\nu$, $G$ be $G_{2\nu}$ or
$G_{2(\nu-1)+2}$, and
$$
O_n({\mathbb F}_q,G)=\{\,T\in{\em \textmd{GL}}_n({\mathbb
F}_q):TG{^tT}\equiv G\,\}.
$$
Let $A=G+{^tG}$ and define $P_{nk}$, $k=0,1,\ldots$, as in $(1)$
and $(2)$. Then
$$
{\mathbb F}_q(X_1,\ldots,X_n)^{O_n({\mathbb F}_q,G)}={\mathbb
F}_q(P_{n0},P_{n1},\ldots,P_{n,n-1}).
$$
\end{Theorem}

\begin{pf}
Let
$$
G_A=\left\{\,T\in\textmd{GL}_n({\mathbb F}_q):TA{^tT}=A\right\}.
$$
Then $G_A\cong Sp_{2\nu}({\mathbb F}_q)$. Let $T\in O_n({\mathbb
F}_q,G)$. Then $TG{^tT}+G$ is alternate, hence,
\begin{eqnarray*}
TA{^tT}+A
&=& T(G+{^tG}){^tT}+(G+{^tG})\\
&=& {^t(TG{^tT}+G)}+(TG{^tT}+G)\\
&=&0.
\end{eqnarray*}
Then $T\in G_A$, thus $O_n({\mathbb F}_q,G)\subseteq G_A$. It
follows that
$$
{\mathbb F}_q(X_1,\ldots,X_n)^{G_A}\subseteq {\mathbb
F}_q(X_1,\ldots,X_n)^{O_n({\mathbb F}_q,G)}.
$$
It is known, as mentioned in the introduction, that ${\mathbb
F}_q(X_1,\ldots,X_n)^{G_A}={\mathbb F}_q(P_{n1},\ldots,P_{nn})$.
Furthermore, note that, for any $T\in \textmd{GL}_n({\mathbb
F}_q)$,
\begin{eqnarray*}
T\in O_n({\mathbb F}_q,G)&\Leftrightarrow&
(X_1,\ldots,X_n)(TG{^tT}+G)\left(
\begin{array}{c}
X_1\\
\vdots\\
X_n
\end{array}
\right)=0\\
&\Leftrightarrow& \sigma_T(P_{n0})=P_{n0}.
\end{eqnarray*}
Hence
\begin{eqnarray*}
{\mathbb F}_q(X_1,\ldots,X_n)^{O_n({\mathbb F}_q,G)}
&=&{\mathbb F}_q(X_1,\ldots,X_n)^{G_A}(P_{n0})\\
&=&{\mathbb F}_q(P_{n0},P_{n1},\ldots,P_{nn}).
\end{eqnarray*}
It remains to show that $P_{nn}\in {\mathbb
F}_q(P_{n0},P_{n1},\ldots,P_{n,n-1})$.

Suppose we work temporarily in ${\mathbb Z}[X_1,\ldots,X_n]$ as we
did before this theorem. Let
$$
\widetilde{P}= \left( \begin{array}{ccccc}
P_{n1}&P_{n2}&\ldots&P_{n,n-1}&P_{nn}\\
2P^{(q)}_{n0}&P^{(q)}_{n1}&\ldots&P^{(q)}_{n,n-2}&P^{(q)}_{n,n-1}\\
P^{(q)}_{n1}&2P^{(q^2)}_{n0}&\ldots&P^{(q^2)}_{n,n-3}&P_{n,n-2}^{(q^2)}\\
&\ldots&\ldots&\ldots&\\
P^{(q)}_{n,n-2}&P_{n,n-3}^{(q^2)}&\ldots&2P_{n0}^{(q^{n-1})}&P^{(q^{n-1})}_{n1}
\end{array} \right),
$$
then $BA{^tB^{(q)}}=\widetilde{P}$. Note that
\begin{eqnarray*}
\det(\widetilde{P})&=& \det\left( \begin{array}{ccccc}
P_{n1}&P_{n2}&\ldots&P_{n,n-1}&0\\
2P^{(q)}_{n0}&P^{(q)}_{n1}&\ldots&P^{(q)}_{n,n-2}&P^{(q)}_{n,n-1}\\
P^{(q)}_{n1}&2P^{(q^2)}_{n0}&\ldots&P^{(q^2)}_{n,n-3}&P_{n,n-2}^{(q^2)}\\
&\ldots&\ldots&\ldots&\\
P^{(q)}_{n,n-2}&P_{n,n-3}^{(q^2)}&\ldots&2P_{n0}^{(q^{n-1})}&P^{(q^{n-1})}_{n1}
\end{array} \right)\\
&&-P_{nn}\det\left( \begin{array}{cccc}
2P^{(q)}_{n0}&P^{(q)}_{n1}&\ldots&P^{(q)}_{n,n-2}\\
P^{(q)}_{n1}&2P^{(q^2)}_{n0}&\ldots&P^{(q^2)}_{n,n-3}\\
&\ldots&\ldots&\\
P^{(q)}_{n,n-2}&P_{n,n-3}^{(q^2)}&\ldots&2P_{n0}^{(q^{n-1})}
\end{array} \right)\\
&=& 2f+\det(P_n)-\det(P_0)P_{nn},
\end{eqnarray*}
where $f \in {\mathbb
Z}[P_{n0}^{(q)},\ldots,P_{n0}^{(q^{n-1})},P_{n1}^{(q^{n-1})},P]$
and
\begin{eqnarray*}
P_n&=&\left( \begin{array}{ccccc}
P_{n1}&P_{n2}&\ldots&P_{n,n-1}&0\\
0&P^{(q)}_{n1}&\ldots&P^{(q)}_{n,n-2}&P^{(q)}_{n,n-1}\\
P^{(q)}_{n1}&0&\ldots&P^{(q^2)}_{n,n-3}&P_{n,n-2}^{(q^2)}\\
&\ldots&\ldots&\ldots&\\
P^{(q)}_{n,n-2}&P_{n,n-3}^{(q^2)}&\ldots&0&P^{(q^{n-1})}_{n1}
\end{array} \right),\\
P_0&=&\left( \begin{array}{cccc}
2P^{(q)}_{n0}&P^{(q)}_{n1}&\ldots&P^{(q)}_{n,n-2}\\
P^{(q)}_{n1}&2P^{(q^2)}_{n0}&\ldots&P^{(q^2)}_{n,n-3}\\
&\ldots&\ldots&\\
P^{(q)}_{n,n-2}&P_{n,n-3}^{(q^2)}&\ldots&2P_{n0}^{(q^{n-1})}
\end{array} \right).
\end{eqnarray*}
Denote the matrix obtained from $P_0$ by replacing all its
diagonal elements by $0$ and deleting its $i$-th row and $i$-th
column simultaneously by $P_{(i)}$. Then $\det(P_0)\equiv 2P_0^*
(\mod 4)$ where
$$
P_0^*=P_{n0}^{(q)}\det(P_{(1)})+P_{n0}^{(q^2)}\det(P_{(2)})+\ldots+P_{n0}^{(q^{n-1})}\det(P_{(n-1)})+l
$$
with $l\in{\mathbb Z}[P]$ such that
$$
2l=\det \left( \begin{array}{cccc}
0&P^{(q)}_{n1}&\ldots&P^{(q)}_{n,n-2}\\
P^{(q)}_{n1}&0&\ldots&P^{(q^2)}_{n,n-3}\\
&\ldots&\ldots&\\
P^{(q)}_{n,n-2}&P_{n,n-3}^{(q^2)}&\ldots&0
\end{array} \right),
$$
(the existence of $l$ follows from \ref{nodd}). Clearly, $P_0^*\in
{\mathbb Z}[P_{n0}^{(q)},\ldots,P_{n0}^{(q^{n-1})},P]$. By
\ref{mainlemma}, $\det(B)\equiv 2S+\Pf(P) (\mod 4)$, hence $\det
(B^{(q)})\equiv 2S^{(q)}+\Pf(P^{(q)}) (\mod 4)$. Then, from
$\det(\widetilde{P})=\det(BA^tB^{(q)})$, we have that
\begin{eqnarray*}
&&
 2f+\det(P_n)-2P_0^*P_{nn}\\
 &\equiv&\det(A)\det(B)\det(B^{(q)}) (\mod 4)\\
 &\equiv&((-1)^{\nu-1}4\alpha^2+(-1)^{\nu})(2S+\Pf(P))(2S^{(q)}+\Pf(P^{(q)})) (\mod 4)\\
 &\equiv&
 (-1)^{\nu}2(S\Pf(P^{(q)})+\Pf(P)S^{(q)})+(-1)^{\nu}\Pf(P)\Pf(P^{(q)})(\mod
 4).
 \end{eqnarray*}
In virtue of \ref{pftimes},
\begin{eqnarray*}
\Pf(P)\Pf(P^{(q)})&=& \det
\left( \begin{array}{ccccc} 0&P_{n1}&P_{n2}&\ldots&P_{n,n-1}\\
-P_{n,n-1}^{(q)}&0&P_{n1}^{(q)}&\ldots&P_{n,n-2}^{(q)}\\
-P_{n,n-2}^{(q^2)}&-P_{n1}^{(q)}&0&\ldots&P_{n,n-3}^{(q^2)}\\
&\ldots&&\ldots&\\
-P_{n1}^{(q^{n-1})}&-P_{n,n-2}^{(q)}&-P_{n,n-3}^{(q^2)}&\ldots&0
\end{array} \right)\\
&=&\det\left( \begin{array}{ccccc}
P_{n1}&P_{n2}&\ldots&P_{n,n-1}&0\\
0&P^{(q)}_{n1}&\ldots&P^{(q)}_{n,n-2}&P^{(q)}_{n,n-1}\\
-P^{(q)}_{n1}&0&\ldots&P^{(q^2)}_{n,n-3}&P_{n,n-2}^{(q^2)}\\
&\ldots&\ldots&\ldots&\\
-P^{(q)}_{n,n-2}&-P_{n,n-3}^{(q^2)}&\ldots&0&P^{(q^{n-1})}_{n1}
\end{array} \right).
\end{eqnarray*}
It is easy to see, by the same argument in the proof of
\ref{nodd}, that
\begin{eqnarray*}
&& \det\left( \begin{array}{ccccc}
P_{n1}&P_{n2}&\ldots&P_{n,n-1}&0\\
0&P^{(q)}_{n1}&\ldots&P^{(q)}_{n,n-2}&P^{(q)}_{n,n-1}\\
-P^{(q)}_{n1}&0&\ldots&P^{(q^2)}_{n,n-3}&P_{n,n-2}^{(q^2)}\\
&\ldots&\ldots&\ldots&\\
-P^{(q)}_{n,n-2}&-P_{n,n-3}^{(q^2)}&\ldots&0&P^{(q^{n-1})}_{n1}
\end{array} \right)-\det(P_n)\\
&\in & 2{\mathbb
Z}[P_{n1}^{(q^{n-1})},P_{n2}^{(q^{n-2})},\ldots,P_{n,n-1}^{(q)},P].
\end{eqnarray*}
Hence there exists $h\in {\mathbb
Z}[P_{n1}^{(q^{n-1})},P_{n2}^{(q^{n-2})},\ldots,P_{n,n-1}^{(q)},P]$
such that\\
$\det(P_n)-(-1)^{\nu}\Pf(P)\Pf(P^{(q)})=2h$. Then
$$
P_0^*P_{nn}\equiv f+(-1)^{\nu-1}(S\Pf(P^{(q)})+\Pf(P)S^{(q)})+h
(\mod 2).
$$

Now pass to ${\mathbb Z}_2[X_1,\ldots,X_n]$, hence, to ${\mathbb
F}_q[X_1,\ldots,X_n]$, thus, $\alpha$ is an element of ${\mathbb
F}_q$. Then we have
\begin{eqnarray*}
P_0^*P_{nn}&=&
f+h+S\Pf(P)^q+\Pf(P)S^q\\
&=& f+h+S(S+\Pf(P))^q+S^q(S+\Pf(P)),
\end{eqnarray*}
where $P_0^*,f,h\in {\mathbb Z}[P_{n0},P_{n1},\ldots,P_{n,n-1}]$.

Note that, as a symmetric polynomial in $S$ and $S+\Pf(P)$,
$S(S+\Pf(P))^q+S^q(S+\Pf(P))$ is a polynomial of
$S+(S+\Pf(P))=\Pf(P)$ and $S(S+\Pf(P))$. But
$S(S+\Pf(P))\in{\mathbb Z}_2[P_{n0},P_{n1},\ldots,P_{n,n-1}]$ by
\ref{mainlemma}. Hence $S(S+\Pf(P))^q+S^q(S+\Pf(P))\in{\mathbb
Z}_2[P_{n0},P_{n1},\ldots,P_{n,n-1}]$.

Finally, let us show that $P_0^*\not=0$. When $n=2$, it is clear
that $P_0^*=P_{20}^q\not=0$. Suppose that $n>2$. We distinguish
the cases $G=G_{2\nu}$ and $G=G_{2(\nu-1)+2}$. When $G=G_{2\nu}$,
set $X_{\nu}=X_n$ in $P_0^*$, then $P_{n0}=P_{n-2,0}+X_n^2$,
$P_{ni}=P_{n-2,i}$, $i=1,2,\ldots$, where $P_{n-2,i}$ is defined
relative to $G_{2(\nu-1)}$ and the indeterminates
$X_1,\ldots,X_{\nu-1},X_{\nu+1},\ldots,X_{2\nu-1}$. Hence
$l,\det(P_{(i)})\in {\mathbb
Z}[X_1,\ldots,X_{\nu-1},X_{\nu+1},\ldots,X_{n-1}]$,
$i=1,\ldots,n-1$. As a polynomial in $X_n$, the leading
coefficient of $P_0^*$ is $\det(P_{(n-1)})$. But
\begin{eqnarray*}
&&\det(P_{(n-1)})\\
&=& \det\left( \begin{array}{cccc}
0&P^q_{n-2,1}&\ldots&P^q_{n-2,n-3}\\
P^q_{n-2,1}&0&\ldots&P^{q^2}_{n-2,n-4}\\
&\ldots&\ldots&\\
P^{q}_{n-2,n-3}&P_{n-2,n-4}^{q^2}&\ldots&0
\end{array} \right)\\
&=& \det\left\{ \left(
\begin{array}{cccc}
X_1^q&X_2^q&\ldots&X_{n-2}^q\\
X_1^{q^2}&X_2^{q^2}&\ldots&X_{n-2}^{q^2}\\
&\ldots&\ldots&\\
X_1^{q^{n-2}}&X_2^{q^{n-2}}&\ldots&X_{n-2}^{q^{n-2}}
\end{array}
\right)\left(
\begin{array}{cc}
0&I^{(\nu-1)}\\
I^{(\nu-1)}&0
\end{array}
\right)\left(
\begin{array}{cccc}
X_1^q&X_1^{q^2}&\ldots&X_1^{q^{n-2}}\\
X_2^q&X_2^{q^2}&\ldots&X_2^{q^{n-2}}\\
&\ldots&\ldots&\\
X_{n-2}^q&X_{n-2}^{q^2}&\ldots&X_{n-2}^{q^{n-2}}
\end{array}
\right)\right\}\\
&=& \det\left(
\begin{array}{cccc}
X_1&X_2&\ldots&X_{n-2}\\
X_1^q&X_2^q&\ldots&X_{n-2}^q\\
&\ldots&\ldots&\\
X_1^{q^{n-3}}&X_2^{q^{n-3}}&\ldots&X_{n-2}^{q^{n-3}}
\end{array}
\right)^{q^2}\\
&\not=&0.
\end{eqnarray*}
Hence $P_0^*\not=0$. Then $P_{nn}\in{\mathbb
F}_q(P_{n0},P_{n1},\ldots,P_{n,n-1})$. When $G=G_{2(\nu-1)+2}$, we
take $G=\left(\begin{array}{cccc}\alpha&1&&\\
&\alpha&&\\
&&0&I^{(\nu-1)}\\
&&&0
\end{array}\right)$ and set $X_{\nu+1}=X_n$, then we also have
$P_{n0}=P_{n-2,0}+X_n^2$, $P_{ni}=P_{n-2,i}$, $i=1,2,\ldots$, but
now $P_{n-2,i}$ is defined relative to $\left(\begin{array}{cccc}\alpha&1&&\\
&\alpha&&\\
&&0&I^{(\nu-2)}\\
&&&0
\end{array}\right)$ and the
indeterminates $X_1,\ldots,X_{\nu},X_{\nu+2},\ldots,X_{2\nu-1}$.
Thus we can draw the same conclusion $P_{nn}\in{\mathbb
F}_q(P_{n0},P_{n1},\ldots,P_{n,n-1})$.
\end{pf}

Now let us come to the case when $n=2\nu+1$ is odd. We need the
following lemma.

\begin{Lemma}
\label{det} Let $n=2\nu$ and ${\mathbb F}_q$ a finite field of
characteristic two. In ${\mathbb F}_q[X_1,\ldots,X_n]$, let
\begin{eqnarray*}
P_{nk}&=& (X_1,\ldots,X_n) \left(
\begin{array}{cc}
0&I^{(\nu)}\\
I^{(\nu)}&0
\end{array}
\right)\left(
\begin{array}{c}
X_1^{q^k}\\
\vdots\\
X_n^{q^k}
\end{array}
\right),\,\,k=1,2,\ldots,\\
D_{i_1i_2\cdots i_n}&=&\det \left(
\begin{array}{cccc}
X_1^{q^{i_1}}&X_1^{q^{i_2}}&\ldots&X_1^{q^{i_n}}\\
X_2^{q^{i_1}}&X_2^{q^{i_2}}&\ldots&X_2^{q^{i_n}}\\
&\ldots&\ldots&\\
X_n^{q^{i_1}}&X_n^{q^{i_2}}&\ldots&X_n^{q^{i_n}}
\end{array}
\right),\,\,0\leq i_1<i_2<\cdots<i_n.
\end{eqnarray*}
Then
$$
D_{i_1i_2\cdots i_n}\in{\mathbb
F}_q[P_{n1},P_{n2},\ldots,P_{n,i_n-i_1}].
$$
\end{Lemma}

\begin{pf}
Since
\begin{eqnarray*}
(D_{i_1i_2\cdots i_n})^2&=&\det\left\{ \left(
\begin{array}{cccc}
X_1^{q^{i_1}}&X_2^{q^{i_1}}&\ldots&X_n^{q^{i_1}}\\
X_1^{q^{i_2}}&X_2^{q^{i_2}}&\ldots&X_n^{q^{i_2}}\\
&\ldots&\ldots&\\
X_1^{q^{i_n}}&X_2^{q^{i_n}}&\ldots&X_n^{q^{i_n}}
\end{array}
\right)\left(
\begin{array}{cc}
0&I^{(\nu)}\\
I^{(\nu)}&0
\end{array}
\right)\left(
\begin{array}{cccc}
X_1^{q^{i_1}}&X_1^{q^{i_2}}&\ldots&X_1^{q^{i_n}}\\
X_2^{q^{i_1}}&X_2^{q^{i_2}}&\ldots&X_2^{q^{i_n}}\\
&\ldots&\ldots&\\
X_n^{q^{i_1}}&X_n^{q^{i_2}}&\ldots&X_n^{q^{i_n}}
\end{array}
\right)\right\}\\
&=& \det\left(
\begin{array}{ccccc}
0&P_{n,i_2-i_1}^{q^{i_1}}&P_{n,i_3-i_1}^{q^{i_1}}&\ldots&P_{n,i_n-i_1}^{q^{i_1}}\\
P_{n,i_2-i_1}^{q^{i_1}}&0&P_{n,i_3-i_2}^{q^{i_2}}&\ldots&P_{n,i_n-i_2}^{q^{i_2}}\\
P_{n,i_3-i_1}^{q^{i_1}}&P_{n,i_3-i_2}^{q^{i_2}}&0&\ldots&P_{n,i_n-i_3}^{q^{i_3}}\\
&\ldots&&\ldots&\\
P_{n,i_n-i_1}^{q^{i_1}}&P_{n,i_n-i_2}^{q^{i_2}}&P_{n,i_n-i_3}^{q^{i_3}}&\ldots&0
\end{array}
\right),
\end{eqnarray*}
we have that
$$
D_{i_1i_2\cdots i_n}=\Pf\left(
\begin{array}{ccccc}
0&P_{n,i_2-i_1}^{q^{i_1}}&P_{n,i_3-i_1}^{q^{i_1}}&\ldots&P_{n,i_n-i_1}^{q^{i_1}}\\
P_{n,i_2-i_1}^{q^{i_1}}&0&P_{n,i_3-i_2}^{q^{i_2}}&\ldots&P_{n,i_n-i_2}^{q^{i_2}}\\
P_{n,i_3-i_1}^{q^{i_1}}&P_{n,i_3-i_2}^{q^{i_2}}&0&\ldots&P_{n,i_n-i_3}^{q^{i_3}}\\
&\ldots&&\ldots&\\
P_{n,i_n-i_1}^{q^{i_1}}&P_{n,i_n-i_2}^{q^{i_2}}&P_{n,i_n-i_3}^{q^{i_3}}&\ldots&0
\end{array}
\right)\in {\mathbb F}_q[P_{n1},P_{n2},\ldots,P_{n,i_n-i_1}].
$$
\end{pf}

\begin{Theorem}
\label{second} Let ${\mathbb F}_q$ be a finite field of
characteristic two, $n=2\nu+1$, $G=G_{2\nu+1}$ and
$$ O_n({\mathbb
F}_q,G)=\{\,T\in{\em \textmd{GL}}_n({\mathbb F}_q):TG{^tT}\equiv
G\,\}.
$$
Let
\begin{eqnarray*}
P_{n-1,0}&=&(X_1,\ldots,X_{n-1})\left(\begin{array}{cc} 0&I^{(\nu)}\\
0&0
\end{array}\right)\left(
\begin{array}{c}
X_1\\
\vdots\\
X_{n-1}
\end{array}
\right)\\
P_{n-1,k}&=&(X_1,\ldots,X_{n-1}) \left(
\begin{array}{cc}
0&I^{(\nu)}\\
I^{(\nu)}&0
\end{array}
\right)\left(
\begin{array}{c}
X_1^{q^k}\\
\vdots\\
X_{n-1}^{q^k}
\end{array}
\right),\,\,k=1,2,\ldots.
\end{eqnarray*}
Then
$$
{\mathbb F}_q(X_1,\ldots,X_n)^{O_n({\mathbb F}_q,G)}={\mathbb
F}_q(P_{n-1,0},P_{n-1,1},\ldots,P_{n-1,n-2},X_n).
$$
\end{Theorem}

\begin{pf} By Dickson's Theorem,
$$
{\mathbb F}_q(X_1,\ldots,X_n)^{\textmd{GL}_n({\mathbb
F}_q)}={\mathbb F}_q(C_{n0},C_{n1},\ldots,C_{n,n-1}),
$$
where $C_{ni}=\frac{D_{ni}}{D_{nn}}$, $i=0,1,\ldots,n-1$, and
$$
D_{ni}=\det\left(
\begin{array}{cccccc}
X_1&X_1^q&\ldots&\widehat{X_1^{q^i}}&\ldots&X_1^{q^n}\\
X_2&X_2^q&\ldots&\widehat{X_2^{q^i}}&\ldots&X_2^{q^n}\\
&\ldots&&&\ldots&\\
X_n&X_n^q&\ldots&\widehat{X_n^{q^i}}&\ldots&X_n^{q^n}
\end{array}
\right),\,\, i=0,1,\ldots,n.
$$
Let
$$
P_{n0}=(X_1,\ldots,X_n)G\left(
\begin{array}{c}
X_1\\
\vdots\\
X_n
\end{array}
\right).
$$
Note that, for any $T\in\textmd{GL}_n({\mathbb F}_q)$,
\begin{eqnarray*}
T\in O_n({\mathbb F}_q,G)&\Leftrightarrow&
(X_1,\ldots,X_n)(TG{^tT}+G)\left(
\begin{array}{c}
X_1\\
\vdots\\
X_n
\end{array}
\right)=0\\
&\Leftrightarrow& \sigma_T(P_{n0})=P_{n0}.
\end{eqnarray*}
Hence
\begin{eqnarray*}
{\mathbb F}_q(X_1,\ldots,X_n)^{O_n({\mathbb F}_q,G)} &=&{\mathbb
F}_q(X_1,\ldots,X_n)^{\textmd{GL}_n({\mathbb
F}_q)}(P_{n0})\\
&=&{\mathbb F}_q(P_{n0},C_{n0},C_{n1},\ldots,C_{n,n-1}).
\end{eqnarray*}
It is easy to see that $P_{n-1,1},\ldots,P_{n-1,n}\in{\mathbb
F}_q(X_1,\ldots,X_n)^{O_n({\mathbb F}_q,G)}$. By \cite[Theorem
7.1]{W}, every $T\in O_n({\mathbb F}_q,G)$ has the form
$$
\left(\begin{array}{ccc} A&B&E\\
C&D&F\\
0&0&1
\end{array}\right),
$$
where $\left(\begin{array}{cc} A&B\\
C&D\end{array}\right)\in Sp_{2\nu}({\mathbb F}_q)$, it follows
that $\sigma_T(X_n)=X_n$ which implies that $X_n\in {\mathbb
F}_q(X_1,\ldots,X_n)^{O_n({\mathbb F}_q,G)}$. Then
$P_{n-1,0}=P_{n0}+X_n^2\in{\mathbb
F}_q(X_1,\ldots,X_n)^{O_n({\mathbb F}_q,G)}$. In virtue of
\ref{det} and
\begin{eqnarray*}
D_{ni}&=&\sum_{j=0}^{i-1}X_n^{q^j}\det \left(
\begin{array}{cccccccc}
X_1&X_1^q&\ldots&\widehat{X_1^{q^j}}&\ldots&\widehat{X_1^{q^i}}&\ldots&X_1^{q^n}\\
X_2&X_2^q&\ldots&\widehat{X_2^{q^j}}&\ldots&\widehat{X_2^{q^i}}&\ldots&X_2^{q^n}\\
&&\ldots&&&\ldots&&\\
X_{n-1}&X_{n-1}^q&\ldots&\widehat{X_{n-1}^{q^j}}&\ldots&\widehat{X_{n-1}^{q^i}}&\ldots&X_{n-1}^{q^n}
\end{array}
\right)\\
&&+\sum_{j=i+1}^nX_n^{q^j}\det \left(
\begin{array}{cccccccc}
X_1&X_1^q&\ldots&\widehat{X_1^{q^i}}&\ldots&\widehat{X_1^{q^j}}&\ldots&X_1^{q^n}\\
X_2&X_2^q&\ldots&\widehat{X_2^{q^i}}&\ldots&\widehat{X_2^{q^j}}&\ldots&X_2^{q^n}\\
&&\ldots&&&\ldots&&\\
X_{n-1}&X_{n-1}^q&\ldots&\widehat{X_{n-1}^{q^i}}&\ldots&\widehat{X_{n-1}^{q^j}}&\ldots&X_{n-1}^{q^n}
\end{array}
\right),\\
&&i=0,1,\ldots,n.
\end{eqnarray*}
we see that, $D_{ni}\in{\mathbb
F}_q[P_{n-1,1},\ldots,P_{n-1,n},X_n]$, $i=0,1,\ldots,n$. But, as
seen in the introduction, $P_{n-1,n}\in{\mathbb
F}_q(P_{n-1,1},\ldots,P_{n-1,n-1})$ and, from the proof of
\ref{mainthm}, we see that $P_{n-1,n-1}\in{\mathbb
F}_q(P_{n-1,0},P_{n-1,1},\ldots,P_{n-1,n-2})$. Hence
$D_{ni}\in{\mathbb
F}_q(P_{n-1,0},P_{n-1,1},\ldots,P_{n-1,n-2},X_n)$, hence,
$C_{ni}\in{\mathbb
F}_q(P_{n-1,0},P_{n-1,1},\ldots,P_{n-1,n-2},X_n)$,
$i=0,1,\ldots,n-1$. Therefore,
$$
{\mathbb F}_q(X_1,\ldots,X_n)^{O_n({\mathbb F}_q,G)}={\mathbb
F}_q(P_{n-1,0},P_{n-1,1},\ldots,P_{n-1,n-2},X_n).
$$
\end{pf}

\medskip
\section{Rational Invariants of Pseudo-Symplectic Groups}
Let ${\mathbb F}_q$ be a finite field of characteristic two. In
this section, we will discuss the rational invariants of
pseudo-symplectic groups.

\begin{Theorem}
\label{s1} Let $n=2\nu+1$, $S=S_{2\nu+1}$ and
$$
P_{n-1,k}=(X_1,\ldots,X_{n-1}) \left(
\begin{array}{cc}
0&I^{(\nu)}\\
I^{(\nu)}&0
\end{array}
\right)\left(
\begin{array}{c}
X_1^{q^k}\\
\vdots\\
X_{n-1}^{q^k}
\end{array}
\right),\,\,k=1,2,\ldots.
$$
Then
$$
{\mathbb F}_q(X_1,\ldots,X_n)^{Ps_n({\mathbb F}_q,S)}={\mathbb
F}_q(P_{n-1,1},\ldots,P_{n-1,n-1},X_n).
$$
\end{Theorem}

\begin{pf}
Let
$$
\overline{P}_{nk}=(X_1,\ldots,X_n)S\left(
\begin{array}{c}
X_1^{q^k}\\
\vdots\\
X_n^{q^k}
\end{array}
\right),\,\,k=0,1,2,\ldots.
$$
Then, as mentioned in the introduction,
$$
{\mathbb F}_q(X_1,\ldots,X_n)^{Ps_n({\mathbb F}_q,S)}={\mathbb
F}_q(\overline{P}_{n0},\overline{P}_{n1},\ldots,\overline{P}_{nn}).
$$
Note that $\overline{P}_{n0}=X_n^2$ and, for any $T\in
Ps_n({\mathbb F}_q,S)$, from
$\sigma_T(\overline{P}_{n0})=\overline{P}_{n0}$, i.e.,
$(\sigma_T(X_n))^2=X_n^2$, we see that $\sigma_T(X_n)=X_n$, thus
$X_n\in {\mathbb F}_q(X_1,\ldots,X_n)^{Ps_n({\mathbb F}_q,S)}$.
Since $\overline{P}_{nk}=P_{n-1,k}+X_n^{q^k+1}$, it follows that
$$
{\mathbb F}_q(X_1,\ldots,X_n)^{Ps_n({\mathbb F}_q,S)}={\mathbb
F}_q(P_{n-1,1},\ldots,P_{n-1,n},X_n).
$$
If we set
$$
A=\left(
\begin{array}{cc}
0&I^{(\nu)}\\
I^{(\nu)}&0
\end{array}
\right),
$$
then, as mentioned in the introduction,
$$
P_{n-1,n}\in {\mathbb F}_q(X_1,\ldots,X_{n-1})^{G_A}={\mathbb
F}_q(P_{n-1,1},\ldots,P_{n-1,n-1}).
$$
Hence $ {\mathbb F}_q(X_1,\ldots,X_n)^{Ps_n({\mathbb
F},S)}={\mathbb F}_q(P_{n-1,1},\ldots,P_{n-1,n-1},X_n)$.
\end{pf}

\begin{Lemma}
\label{pfaffian} Let $f_{ij}\in {\mathbb F}_q[X_1,\ldots,X_n]$,
$1\leq i\leq j\leq m$. Then there exists $f\in {\mathbb
F}_q[f_{ij},1\leq i\leq j\leq m]$ such that
$$
\det \left( \begin{array}{cccc} f^2_{11}&f_{12}&\ldots&f_{1m}\\
f_{12}&f^2_{22}&\ldots&f_{2m}\\
&\ldots&\ldots&\\
f_{1m}&f_{2m}&\ldots&f^2_{mm} \end{array} \right)=f^2.
$$
\end{Lemma}

\begin{pf} Let
$$
B=\left( \begin{array}{cccc} f^2_{11}&f_{12}&\ldots&f_{1m}\\
f_{12}&f^2_{22}&\ldots&f_{2m}\\
&\ldots&\ldots&\\
f_{1m}&f_{2m}&\ldots&f^2_{mm} \end{array} \right).
$$
Let $r$ be the number of nonzero elements in the diagonal of $B$.
We prove the lemma by induction on $r$. If $r=0$, then $B$ is
alternate, hence $\det(B)=(\Pf(B))^2$. Now assume that $r>0$. We
may assume that $f_{11}\not=0$. Then
$\det(B)=f^2_{11}\det(B_1)+\det(B_2)$, where
$$
B_1=\left( \begin{array}{cccc} f^2_{22}&f_{23}&\ldots&f_{2m}\\
f_{23}&f^2_{33}&\ldots&f_{3m}\\
&\ldots&\ldots&\\
f_{2m}&f_{3m}&\ldots&f^2_{mm} \end{array} \right),\,\,
B_2=\left( \begin{array}{cccc} 0& f_{12}&\ldots&f_{1m}\\
f_{12}&f^2_{22}&\ldots&f_{2m}\\
&\ldots&\ldots&\\
f_{1m}&f_{2m}&\ldots&f^2_{mm} \end{array} \right).
$$
By induction assumption, $\det(B_1)=f_1^2$, $\det(B_2)=f_2^2$, for
some $f_1,f_2\in {\mathbb F}_q[f_{ij},1\leq i\leq j\leq m]$. Let
$f=f_{11}f_1+f_2$. Then $\det(B)=f^2$, as required.
\end{pf}

\begin{Theorem}
\label{s2} Let $n=2\nu+2$, $S=S_{2\nu+2}$ and
$$
P_{nk}=(X_1,\ldots,X_n)S\left(
\begin{array}{c}
X_1^{q^k}\\
\vdots\\
X_n^{q^k}
\end{array}
\right),\,\,k=0,1,2,\ldots.
$$
Then
$$
{\mathbb F}_q(X_1,\ldots,X_n)^{Ps_n({\mathbb F}_q,S)}={\mathbb
F}_q(P_{n1},\ldots,P_{n,n-1},X_n).
$$
\end{Theorem}

\begin{pf} As in the proof of \ref{s1}, we have
$$
{\mathbb F}_q(X_1,\ldots,X_n)^{Ps_n({\mathbb F}_q,S)}={\mathbb
F}_q(P_{n1},\ldots,P_{nn},X_n).
$$
Let
$$
D_n=\det\left(
\begin{array}{cccc}
X_1&X_2&\ldots&X_n\\
X_1^q&X_2^q&\ldots&X_n^q\\
&\ldots&\ldots&\\
X_1^{q^{n-1}}&X_2^{q^{n-1}}&\ldots&X_n^{q^{n-1}}
\end{array}
\right).
$$
Then
\begin{eqnarray*}
D_n^2&=& \det\left\{ \left(
\begin{array}{cccc}
X_1&X_2&\ldots&X_n\\
X_1^q&X_2^q&\ldots&X_n^q\\
&\ldots&\ldots&\\
X_1^{q^{n-1}}&X_2^{q^{n-1}}&\ldots&X_n^{q^{n-1}}
\end{array}
\right)S\left(
\begin{array}{cccc}
X_1&X_1^q&\ldots&X_1^{q^{n-1}}\\
X_2&X_2^q&\ldots&X_2^{q^{n-1}}\\
\vdots&\vdots&&\vdots\\
X_n&X_n^q&\ldots&X_n^{q^{n-1}}
\end{array}
\right)\right\}\\
&=& \det\left(
\begin{array}{cccc}
P_{n0}&P_{n1}&\ldots&P_{n,n-1}\\
P_{n1}&P_{n0}^q&\ldots&P_{n,n-2}^q\\
&\ldots&\ldots&\\
P_{n,n-1}&P_{n,n-2}^q&\ldots&P_{n0}^{q^{n-1}}
\end{array}
\right)
\end{eqnarray*}
and $P_{n0}=X_n^2$. It follows from \ref{pfaffian} that
$D^2_n=f^2$ for some $f\in {\mathbb
F}_q[P_{n1},\ldots,P_{n,n-1},X_n]$. Hence $D_n\in {\mathbb
F}_q[P_{n1},\ldots,P_{n,n-1},X_n]$. On the other hand, we have
\begin{eqnarray*}
D_n^{q+1}&=& \det\left\{ \left(
\begin{array}{cccc}
X_1^q&X_2^q&\ldots&X_n^q\\
X_1^{q^2}&X_2^{q^2}&\ldots&X_n^{q^2}\\
&\ldots&\ldots&\\
X_1^{q^{n}}&X_2^{q^{n}}&\ldots&X_n^{q^{n}}
\end{array}
\right)S\left(
\begin{array}{cccc}
X_1&X_1^q&\ldots&X_1^{q^{n-1}}\\
X_2&X_2^q&\ldots&X_2^{q^{n-1}}\\
\vdots&\vdots&&\vdots\\
X_n&X_n^q&\ldots&X_n^{q^{n-1}}
\end{array}
\right)\right\}\\
&=& \det\left(
\begin{array}{ccccc}
P_{n1}&P_{n0}^q&P_{n1}^q&\ldots&P_{n,n-2}^q\\
P_{n2}&P_{n1}^q&P_{n0}^{q^2}&\ldots&P_{n,n-3}^{q^2}\\
&\ldots&\ldots&\ldots&\\
P_{n,n-1}&P_{n,n-2}^q&P_{n,n-3}^{q^2}&\ldots&P_{n0}^{q^{n-1}}\\
P_{nn}&P_{n,n-1}^q&P_{n,n-2}^{q^2}&\ldots&P_{n1}^{q^{n-1}}
\end{array}
\right).
\end{eqnarray*}
Let
$$
K=\det\left(
\begin{array}{cccc}
P_{n0}&P_{n1}&\ldots&P_{n,n-2}\\
P_{n1}&P_{n0}^{q}&\ldots&P_{n,n-3}^{q}\\
&\ldots&\ldots&\\
P_{n,n-2}&P_{n,n-3}^{q}&\ldots&P_{n0}^{q^{n-2}}
\end{array}
\right).
$$
Then $D_n^{q+1}=K^qP_{nn}+g$, for some $g\in {\mathbb
F}_q[P_{n1},\ldots,P_{n,n-1}]$. Note that, $K\in {\mathbb
F}_q[P_{n1},\ldots,P_{n,n-2},X_n]$. Furthermore, since
\begin{eqnarray*}
&&K|_{X_{n-1}=X_n}\\
&=& \det \left\{ \left(
\begin{array}{ccccc}
X_1&\ldots&X_{n-2}&X_{n-1}&X_{n-1}\\
X_1^q&\ldots&X_{n-2}^q&X_{n-1}^q&X_{n-1}^q\\
&\ldots&&\ldots&\\
X_1^{q^{n-2}}&\ldots&X_{n-2}^{q^{n-2}}&X_{n-1}^{q^{n-2}}&X_{n-1}^{q^{n-2}}
\end{array}
\right)S\left(
\begin{array}{cccc}
X_1&X_1^q&\ldots&X_1^{q^{n-2}}\\
\vdots&\vdots&&\vdots\\
X_{n-2}&X_{n-2}^q&\ldots&X_{n-2}^{q^{n-2}}\\
X_{n-1}&X_{n-1}^q&\ldots&X_{n-1}^{q^{n-2}}\\
X_{n-1}&X_{n-1}^q&\ldots&X_{n-1}^{q^{n-2}}
\end{array}
\right)\right\}\\
&=& \det \left\{ \left(
\begin{array}{ccccc}
X_1&\ldots&X_{n-2}&X_{n-1}&0\\
X_1^q&\ldots&X_{n-2}^q&X_{n-1}^q&0\\
&\ldots&&\ldots&\\
X_1^{q^{n-2}}&\ldots&X_{n-2}^{q^{n-2}}&X_{n-1}^{q^{n-2}}&0
\end{array}
\right)\left(
\begin{array}{cccc}
0&I^{(\nu)}&&\\
I^{(\nu)}&0&&\\
&&1&\\
&&&1
\end{array}
\right)\left(
\begin{array}{cccc}
X_1&X_1^q&\ldots&X_1^{q^{n-2}}\\
\vdots&\vdots&&\vdots\\
X_{n-2}&X_{n-2}^q&\ldots&X_{n-2}^{q^{n-2}}\\
X_{n-1}&X_{n-1}^q&\ldots&X_{n-1}^{q^{n-2}}\\
X_{n-1}&X_{n-1}^q&\ldots&X_{n-1}^{q^{n-2}}
\end{array}
\right)\right\}\\
&=& \det \left\{ \left(
\begin{array}{cccc}
X_1&\ldots&X_{n-2}&X_{n-1}\\
X_1^q&\ldots&X_{n-2}^q&X_{n-1}^q\\
&\ldots&\ldots&\\
X_1^{q^{n-2}}&\ldots&X_{n-2}^{q^{n-2}}&X_{n-1}^{q^{n-2}}
\end{array}
\right)\left(
\begin{array}{ccc}
0&I^{(\nu)}&\\
I^{(\nu)}&0&\\
&&1
\end{array}
\right)\left(
\begin{array}{cccc}
X_1&X_1^q&\ldots&X_1^{q^{n-2}}\\
\vdots&\vdots&&\vdots\\
X_{n-2}&X_{n-2}^q&\ldots&X_{n-2}^{q^{n-2}}\\
X_{n-1}&X_{n-1}^q&\ldots&X_{n-1}^{q^{n-2}}
\end{array}
\right)\right\}\\
&=& \det\left(
\begin{array}{cccc}
X_1&X_2&\ldots&X_{n-1}\\
X_1^q&X_2^q&\ldots&X_{n-1}^q\\
&\ldots&\ldots&\\
X_1^{q^{n-2}}&X_2^{q^{n-2}}&\ldots&X_{n-1}^{q^{n-2}}
\end{array}
\right)^2\\
&\not=&0,
\end{eqnarray*}
we see that $K\not=0$, thus $P_{nn}\in {\mathbb
F}_q(P_{n1},\ldots,P_{n,n-1},X_n)$. Hence
$$
{\mathbb F}_q(X_1,\ldots,X_n)^{Ps_n({\mathbb F}_q,S)}={\mathbb
F}_q(P_{n1},\ldots,P_{n,n-1},X_n).
$$
\end{pf}

\medskip
\medskip
\noindent
Zhongming Tang\\
Department of Mathematics\\
Suzhou University\\
Suzhou 215006\\
P.\ R.\ China\\
E-mail: zmtang@@suda.edu.cn

\medskip
\noindent
Zhe-xian Wan\\
Academy of Mathematics and System Sciences\\
Chinese Academy of Science\\
Beijing 100080\\
P.\ R.\ China\\
E-mail: wan@@amss.ac.cn

\end{document}